\input AHTOH-E.STY
\hfuzz 3.6pt
\UDC{
512.543.72 
+
512.544.33 
}

\title{
Virtually free finite-normal-subgroup-free groups are strongly verbally 
closed
}

\author{%
Anton A. Klyachko
\quad
Andrey M. Mazhuga
\quad
Veronika Yu. Miroshnichenko
}
\address{
\myAddressW
\quad
mazhuga.andrew@yandex.ru
\quad
werunik179@gmail.com
}

\grantsFirstSecond{\RFBR15-01-05823}

\abstract{%
Any
virtually free group $H$
containing no non-trivial finite normal subgroup
(e.g.,
the infinite dihedral group) is a retract of any
finitely generated group
containing $H$ as a verbally closed subgroup.
}

\s 0.
Introduction

A subgroup $H$ of a group $G$ is called \emph{verbally closed}
[MR14]
(%
see also
[Rom12],
[RKh13],
[Mazh17],
\hbox{[KlMa18]},
[Mazh18])
if
any equation of the form
$$
w(x_1,x_2,\dots)=h,
\qbox{where $w$ is an element of the free group
$F(x_1,x_2,\dots)$ and $h\in H$,
}
$$
having a solution in $G$ has a solution in $H$.
If each finite system of equations with coefficients from~$H$
$$
\{w_1(x_1,x_2,\dots)=1, \dots, w_m(x_1,x_2,\dots)=1\},
\qbox{where $w_i\in H*F(x_1,x_2,\dots)$,}
$$
having a solution in $G$ has a solution in $H$, then the subgroup $H$
is
called \emph{algebraically closed} in $G$.

Surely, the algebraic closedness is stronger than the verbal
closedness.
However, these properties turn out to be
equivalent in many cases. A group $H$ is called
\emph{strongly verbally closed}
[Mazh18]
if it is algebraically closed
in any group containing $H$ as a verbally closed
subgroup. (Thus, verbal closedness is a property of
a subgroup, while the strong verbal closedness is a property of an abstract
group.)
For example, the following groups are strongly verbally closed:
\-
all abelian groups [Mazh18];
\-
all free groups [KlMa18];
\-
the fundamental groups of all connected surfaces, except possibly
the Klein bottle [Mazh18].



\enditem
The main result of this paper can be stated as follows.

\Theorem 1.
The following groups are strongly verbally closed:
\item{\rm1)}
all 
virtually free group containing no
non-trivial finite normal subgroups;
\item{\rm2)}
all free products
$\zvezda\limits_{i\in I}\!H_i$,
where the set $I$ is finite or infinite,
$|I|>1$, and $H_i$ are
nontrivial groups satisfying
nontrivial laws.

\enditem
A large part of Theorem 1 was known earlier:
in [Mazh18], Assertion 2)
was proved
for all non-dihedral groups
under the additional condition that $I$ is finite;
in [KlMa18],
the strong verbal closedness was proved for
 all infinite virtually free non-dihedral groups containing
no infinite abelian noncyclic subgroups.
Paper [KlMa18]
contains 
also
examples of 
virtually free groups that are not
strongly verbally closed. 

Actually, most of this paper is
devoted to the proof of the 
following particular case of (both assertions of) 
Theorem~1:
\disp{\sl
the infinite dihedral group is strongly verbally closed.%
}
For all non-dihedral groups, Theorem 1 is relatively easily derived 
from known facts (in Section 1).

The difficulty with the infinite dihedral group is that it is
``too abelian" to apply
sophisticated tools based on the Lee words [Lee02] (see
[MR14],
[Mazh17],
[KlMa18],
[Mazh18]);
on the other hand, it is ``too nonabelian" to apply
simple arguments (see
[KlMa18],
[Mazh18]).
Of course, the dihedral group is metabelian and this is the basis of our
approach. In Section 3, we give an ``explicit" criterion 
for an infinite dihedral subgroup
to be verbally
(and algebraically) closed.
This criterion is
similar (in
some sense) to
the following simple fact
about abelian subgroups, which can be easily derived
from a result proved in [Mazh18]:
\disp{\sl
an abelian subgroup $H$ is verbally (and algebraically) closed in
a group $G$
if and only if its intersection with the commutator subgroup $G'$
of $G$ is trivial and the image $H$ in the quotient group $G/G'$ is
pure (servant).
}

The notion of algebraic closedness can be
characterised in structural language
if the group $H$ is \emph{equationally Noetherian}, i.e. any system
of
equations over $H$ with finitely many unknowns is equivalent to its
finite subsystem.
Namely,
the algebraic closedness in this case is equivalent to the
``local retractness" (see Section 1):

\disp{
\sl
An equationally Noetherian
subgroup $H$ of a
group $G$ is algebraically closed
in $G$ if and only if it is a retract of
each
finitely generated over $H$
subgroup of $G$ \rm(i.e. a subgroup of the form $\gp{H,X}$,
where $X\subseteq G$ is a finite set).
}%
All virtually free groups (including the infinite dihedral group)
are
equationally Noetherian [KlMa18]. Therefore, Theorem 1 implies that
\disp{\sl
each 
virtually free group $H$
containing no
finite non-trivial normal subgroup
is a retract of every finitely generated group containing $H$
as a verbally closed subgroup.
}

In Section 1, we prove some auxiliary facts
that allow us to reduce the proof of the main theorem
to the case, where the group $G$ containing a verbally closed
dihedral subgroup $H$ is an extension
of an abelian
group $Q$ by an elementary abelian 2-group $C$; and so
$Q$
is a $C$-module.
Section 2 contains general information about such modules.
In Section 3, we state and prove 
a criterion for 
an infinite dihedral subgroup to be 
algebraically (and verbally) closed.

In the last section, we consider an example
illustrating the
main step of the proof.
This is essentially the simplest example of the situation, where
algebraic unclosedness is almost obvious, while the proof of
verbal unclosedness requires a nontrivial argument. We tried to make
the last section independent; so, readers may read
this section first. We
warn these readers however that the example
considered there
shows some, but not all,
difficulties we face.

\smallskip
\noindent
{\bf Notation},
which we use, is mainly standard. Note only that,
if $X$ is a subset of a group, then $|X|$,
$\gp X$,
and $C(X)$
is the cardinality of $X$, the
subgroup generated by $X$,
and
the centraliser of~$X$,
respectively.
The letters $\Z$ and $\Q$ denote the sets of integers and rationals.
The cyclic group of order $k$ generated by an element $x$ is denoted
$\gp x_k$.
The free group with a basis $x_1,\dots,x_n$ is denoted by
$F(x_1,\dots,x_n)$ (or $F_n$).

\s 1.
Auxiliary lemmata

\Proposition 1.
An equationally Noetherian subgroup $H$ is algebraically closed
in a group $G$ if and only if $H$ is a retract of each
finitely generated over $H$
subgroup of $G$.

\Proof
Suppose that $H$ is algebraically closed in $G$. Then
$H$ is algebraically closed in any subgroup $\widetilde{G}$ of $G$
containing~$H$. If $\widetilde{G}$ is finitely generated over $H$,
then
$H$ is a retract of $\widetilde{G}$ by virtue of the following fact
[MR14]:
\disp{%
\hfuzz15pt
\narrower
\narrower
\narrower
\narrower
\narrower
\sl
An equationally Noetherian algebraically closed subgroup~$H$
of a finitely generated over $H$ group is a retract.
}
Now, suppose that $H$ is a retract of each finitely generated over $H$
subgroup
of $G$ and a system of equations
$
S=\Bigl\{w_1(x_1,\dots,x_n)=1,\dots,w_m(x_1,\dots,x_n)=1\Bigr\},
$
where $w_i\in F(x_1,\dots,x_n)\ast H$, has a solution
$g_1,\dots,g_n$ in $G$.
There exists a retraction
$\rho\:{}\gp{H,g_1,\dots,g_n}\to H$. Therefore,
$\rho(g_1),\dots,\rho(g_n)$ is a solution to~$S$ in $H$
as required.

\Proposition 2.
If each finite subset of a group $H$ is contained
in a strongly verbally closed subgroup which is also verbally closed in
$H$, then $H$ is strongly verbally closed.

\Proof
Suppose that the group $H$ is
verbally closed
in an
overgroup $G$, and
some finite system of equations $S$ with coefficients from~$H$ has
a solution in $G$. The set of coefficients of $S$ is contained in
some subgroup $H_1\subseteq H$ which is strongly verbally
closed and verbally closed in $H$. The latter means that $H_1$ is verbally
closed in $G$ (because verbal closedness is a transitive property).
Now, strong verbal closedness of $H_1$ implies solvability of the system
$S$ in $H_1\subseteq H$ as required.

\Corollary.
Assertion {\rm2)} of Theorem 1 holds for all
non-dihedral
groups.

\Proof
Proposition 2 allows us to reduce the proof to the case, where the set
$I$ is finite, because
any subproduct $H_1=\zvezda\limits_{i\in J}\!H_i$ is verbally
closed
(and even a retract)
in the group~$H=\zvezda\limits_{i\in I}\!H_i$ for any $J\subseteq I$.

For finite $I$, Assertion 2) of Theorem 1
was proved for all non-dihedral groups in [Mazh18].


\Proposition 3.
Assertion {\rm1)} of Theorem 1 holds for all non-dihedral groups.

\Proof
If a group $H$ is virtually cyclic, then it contains a finite
normal subgroup $K$ such that the quotient $H/K$ is either 
trivial,
infinite
cyclic,
or 
infinite
dihedral [Sta71]. 
The normal finite subgroup $K$ must be trivial by the 
condition of Theorem 1.
All abelian groups are strongly verbally closed [Mazh18]
and we are done in this case.

If $H$ is not virtually cyclic, then it contains a
non-abelian free normal subgroup $F$ of a finite index $n$.
This means that the centraliser of the verbal subgroup
$V=\gp{\{h^n\;|\;h\in H\}}\subseteq F$ is normal and finite
(because the centraliser of any nonabelian free subgroup
in a virtually free group is finite).
Hence, $C(V)=\1$ by the condition of Theorem 1.
Therefore, the centraliser of some finitely generated subgroup
$V_1\subseteq V$ is trivial (because again the centraliser of a 
nonabelian
free 
subgroup in a virtually free group is finite 
and a descending chain of finite subgroups stabilises). 
It remains to apply 
the following fact ([Mazh18], Corollary 1):  
\disp{\sl 
if a verbal 
subgroup $V$ of a group $H$ is free nonabelian and the centraliser of a 
finitely generated subgroup of $V$ is trivial, then $H$ is strongly 
verbally closed.
}

The rest of this paper is devoted to the proof that 
the infinite dihedral group is strongly verbal
closed.

\Lemma 1 {\rm([RKh13], Lemma 1.1)}.
If $V(G)$ is a verbal subgroup of a group $G$, and $H$ is a verbally
closed subgroup of $G$, then the image of $H$ under the natural
homomorphism $G\to G/V(G)$ is verbally closed.

\s2.
Commuting involutions on
abelian groups

Suppose that a finite elementary abelian 2-group $C$
(a finite direct power of the two-element group)
acts by automorphisms on a finitely generated
abelian group,
i.e. $Q$ is a $C$-module.
Let $X$ be the set of all homomorphisms (\emph{characters})
$\chi\:C\to\{\pm1\}$ and
$$
Q_\chi=\{q\in Q\;|\; cq={\chi(c)q} \hbox{ for all } c\in C\}.
$$
There is a natural homomorphism
from $Q$ to the additive group of the vector space
(over $\Q$)
$\Q\tensor Q$ sending $q\in Q$ to $1\tensor q$ (the tensor product
is over $\Z$).
The kernel of this homomorphism is the torsion part $T(Q)$ of $Q$.
The action of $C$ on $Q$ extends naturally to a linear
representation $C\to\GL(\Q\tensor Q)$.
This representation is completely reducible and the irreducible
representations are
one-dimensional (the characters of $C$). Thus,
$\Q\tensor Q=\bigoplus\limits_{\chi\in X}(\Q\tensor Q_\chi)$.
The natural projection~$\Q\tensor Q\to\Q\tensor Q_\chi$
is denoted by $p_\chi$.
The vectors $p_\chi(v)\in\Q\tensor Q_\chi$ are called
\emph{$\chi$-components} of the vector~$v\in\Q\tensor Q$ and denoted
by $v_\chi$.
Clearly $\bigoplus\limits_{\chi\in X}(1\tensor
Q_\chi)\subseteq 1\tensor Q\subseteq
\bigoplus\limits_{\chi\in X}p_\chi(1\tensor Q)$
and, if one of these inclusions is an equality, then
the other is an equality; in this case, we say that the $C$-module
$1\tensor Q$ is \emph{decomposable}. In the general case,
the $C$-module~$\bigoplus\limits_{\chi\in X}p_\chi(1\tensor Q)$ can be
called the \emph{decomposable closure} of the module $1\tensor Q$.

We need the following simple formula
valid for any
character $\chi$ and any
$q\in Q$:
$$
1\tensor\(\(\prod_{c\in C}\Bigl(1+\chi(c)c\Bigr)\)\cdot q\)=
\(2^{|C|}\tensor q\)_\chi
\qqbox{and}
|T(Q)|\cdot\prod_{c\in C}\Bigl(1+\chi(c)c\Bigr)\cdot q=
\(2^{|C|}\cdot|T(Q)|\cdot q\)_\chi,
\eqno{(1)}
$$
(where $x_\chi\in Q_\chi$ in the right-hand side 
of the latter equality
are 
components of an element $x\in\bigoplus Q_\chi\subseteq Q$).
Here, the second equality follows from the first one,
because the kernel of the homomorphism
$q\mapsto1\tensor q$ is the torsion part of $Q$
(and the first equality shows that the module
$2^{|C|}\cdot|T(Q)|\cdot Q$ is contained in the direct sum
of $\chi$-components of $Q$).
To prove the first equality note  that the
$\chi$-component of the element $1\tensor q$
in the left-hand side is
multiplied by
two $|C|$ times. As for the other components, they vanish,
because, for each character $\chi'\ne\chi$, there exists
$c\in C$ such that $\chi(c)=-\chi'(c)$. Formula (1) implies
that, for all $q\in Q$,
we have the equality
$$
1\tensor\(\sum_{\chi\in X}\(\prod_{c\in C}\Bigl(1+\chi(c)c\Bigr)\)\cdot q\)
=
2^{|C|}\tensor q
\qqbox{and}
|T(Q)|\cdot
\sum_{\chi\in X}\(\prod_{c\in C}\Bigl(1+\chi(c)c\Bigr)\)\cdot q
=
|T(Q)|\cdot2^{|C|}\cdot q.
\eqno{(2)}
$$
We also need the following simple generalisation of formula (1):
{\sl
If $\phi\:C\to\^C$ is an epimorphism from one finite elementary
abelian 2-group onto another
and $\^q\in\^Q$ is an element of a decomposable $\^C$-module $\^Q$, then,
for any
character $\chi$ of $C$,
$$
1\tensor\(\(\prod_{c\in C}\Bigl(1+\chi(c)\phi(c)\Bigr)\)\cdot\^q\)
=
\cases{
2^{|C|}\tensor\^q_{\^\chi},
& if $\chi=\^\chi\circ\phi$;
\cr\cr
0,
& if
$\chi\ne\^\chi\circ\phi$ for any
character $\^\chi\:{\^C}\to\{\pm1\}$.
}
\eqno{(*)}
$$
}

We call an element $q\in Q$ \emph{simple} if, for some character
$\chi$,
the $\chi$-component
$(1\tensor q)_\chi=p_\chi(1\tensor q)$ is a primitive
element of the free abelian group $p_\chi(1\tensor Q)$, i.e.
$p_\chi(1\tensor q)\notin k\cdot p_\chi(1\tensor Q)=p_\chi(k\tensor Q)$
for $k\in\Z\setminus\{\pm1\}$.

\proclaim{Simple-element Lemma}.
An element $q\in Q$ is simple if and only if
its order is infinite and
the group $Q$ decomposes into a direct sum
$Q=\gp q\oplus M$, where the subgroup $M\subset Q$ is a $C$-submodule,
i.e. $cm\in M$ for all $c\in C$ and $m\in M$.

\Proof
Let $q$ be a simple element, i.e.
$(1\tensor q)_\chi$ is a primitive element of the group
$p_\chi(1\tensor Q)$.
Clearly, this implies that the order of $q$ is infinite
(otherwise $1\tensor q=0$).
Moreover,
$p_\chi(1\tensor Q)=\gp{(1\tensor q)_\chi}\oplus D$,
for some subgroup $D$ (because a primitive element
of a free abelian group is contained in some basis). This
means that
$
1\tensor Q=
\gp{1\tensor q}\oplus
\bigl(p_\chi^{-1}(D)\cap (1\tensor Q)\bigr).
$
Moreover, the group
$A=p_\chi^{-1}(D)$ is a
$C$-submodule. Therefore,
$Q=\gp q\oplus\psi^{-1}\(A\)$, where $\psi\:Q\to\Q\tensor Q$
is the natural homomorphism sending $x$ to $1\tensor x$.

The proof of the other direction is left to readers as
an exercise (we are not going to use it).


\medskip

\Example.
Let the group $C=\gp c_2$ act on $Q=\Z\oplus\Z$
by
permutations of
coordinates. There are two characters: $X=\{\chi_+,\chi_-\}$, where
$\chi_+(c)=1$ and $\chi_-(c)=-1$. Moreover,
$Q_{\chi_+}=\gp{(1,1)}$ and $Q_{\chi_-}=\gp{(1,-1)}$.
The subgroup~$Q_{\chi_+}\oplus Q_{\chi_-}$ has index two in $Q$.
The decomposable closure is
$$
\=Q=\gp{\({1\over2},{1\over2}\)}\oplus\gp{\({1\over2},-{1\over2}\)}=
\{(x,y)\in\Q^2\;|\;x+y\in\Z\ni x-y\}.
$$
The projections $p_{\chi_\pm}\:\=Q\to\(\=Q\)_{\chi_\pm}$
are $p_{\chi_+}(x,y)\mapsto\({x+y\over2},{x+y\over2}\)$
and $p_{\chi_-}(x,y)\mapsto\({x-y\over2},{y-x\over2}\)$.
The element
$(2,5)$ is not simple, because
$(2,5)=7\cdot\({1\over2},{1\over2}\)-3\cdot\({1\over2},-{1\over2}\)$.
Actually, it is easy to
show that an element~$(x,y)\in Q$ is simple if and only if $x\pm y=\pm1$
for some choice of signs.

\s 3.
Algebraically closed infinite dihedral subgroups

Consider a finitely generated group $G$ whose
subgroup $Q=\gp{\{g^2\;|\;g\in G\}}$
generated by squares of all elements is
abelian. The finite elementary abelian 2-group $C=G/Q$ acts on $Q$
by automorphisms
$$
(gQ)\circ q\:=gqg^{-1}
\qbox{(this is well-defined, because $Q$ is abelian)}.
$$
So, $Q$ is a $C$-module and we can apply
the results of the previous section. Note only that now we
stick to multiplicative notation, i.e. we write, e.g.,
$cq_1c^{-1}q_2^2$ instead of $cq_1+2q_2$.
To simplify formulae, we set $\~q\:=q^{|T(Q)|}$.
Formula (1) takes the form
$$
w_\chi(\~q)\:=
f_\chi(c_1,f_\chi(c_2,\dots,f_\chi(c_{|C|},\~q)\dots))
=
\((\~q)^{2^{|C|}}\)_\chi,
\qbox{where $\left\{c_1,\dots,c_{|C|}\right\}=C$}
\eqno{(1')}
$$
and
$f_\chi(gQ,x)\:=xgx^{\chi(gQ)}g^{-1}$ is a ``skew commutator"
(which is well defined, i.e. does not depend on
the choice of the representative $g$ of the coset $c=gQ$).

An analogue of formula (2) takes the form
$$
\prod_{\chi\in X}w_\chi(\~q)=(\~q)^{2^{|C|}}
\qbox{for all $q\in Q$}.
\eqno{(2')}
$$

The strong verbal closedness of the infinite dihedral group
follows immediately from Proposition 1 and the following theorem.

\Theorem 2.
If $H=\gp b_2\semitimes \gp a_\infty$ is an
infinite dihedral subgroup of
a finitely generated
group~$G$,
then the following conditions are equivalent:
\item{\rm1)}
$H$ is verbally closed in $G$;
\item{\rm2)}
$H$ is algebraically closed in $G$;
\item{\rm3)}
$H$ is a retract of $G$ {\rm(i.e. the image
of an endomorphism $\rho$
such that $\rho\circ\rho=\rho$)};
\item{\rm4)}
$a^2Q'$ is a simple element of the $G/Q$-module $Q/Q'$,
where $Q=\gp{\{g^2\;|\;g\in G\}}$.

\Proof
To prove the implication $4)\imp3)$,
first
note that $H\cap Q'=\1$ (otherwise~$a^2Q'$
is of finite order in $Q/Q'$ and is not simple).

Now Simple-element Lemma
provides us with a normal in $G/Q'$ subgroup $M\subset Q/Q'$
such that $Q/Q'=\gp{a^2Q'}\times M$.
The composition of the natural homomorphisms
$G\to G/Q'\to (G/Q')/M=G_1$ is injective on $H$;
the obtained group~$G_1$ is virtually cyclic
(the subgroup generated by squares of all its elements is generated by
element $a^2$).
It is well known that any virtually cyclic group contains a finite
normal subgroup $N$ such that the quotient group is either cyclic or
dihedral (see, e.g., [Sta71]).
Therefore, the composition of homomorphisms
$G\to G/Q'\to (G/Q')/M=G_1\to G_1/N$ is the required retraction onto $H$
(here, we use that the infinite dihedral group has no
finite normal subgroups; and, in $G_1$,
the subgroup generated by squares of all its elements is generated by $a^2$).

The implications
$3)\imp2)\imp1)$ are general facts valid for any groups (see
Introduction).

\medskip

It remains to prove the implication $1)\imp4)$.
By Lemma 1, we can assume
that $G$ satisfies the law $[x^2,y^2]=1$, i.e.
the subgroup $Q$ generated by squares of all elements of $G$ is
abelian (and finitely generated, because it has finite
index in a finitely generated group~$G$).
Indeed, taking the quotient group of $G$
by the commutator subgroup of the subgroup generated by squares of all
elements affects neither $H$, nor condition 1),
(by Lemma 1), nor condition 4).

Suppose that the element $a^2$ is not simple. This means that
$$
\~{(a^2)}\:=a^{2|T(Q)|}=\prod_{\chi\in X}\~{q(\chi)}_\chi^{k_\chi},
\qbox{for some
$k_\chi\in\Z\setminus\{\pm1\}$
and
$q(\chi)\in Q$}
$$
(where $x_\chi$ are components of 
$x\in |T(Q)|Q\subseteq \Q\tensor Q=
\bigoplus\limits_\chi(\Q\tensor Q_\chi)$).
Formula~$(2')$ gives the equality
$$
\prod_{\chi\in X}
w_\chi\Bigl(\~{(q(\chi)})^{k_\chi}\Bigr)
=
\~{\(a^2\)}^{2^{|C|}}.
\eqno{(3)}
$$

Now, we decompose the finite elementary abelian 2-group $G/Q$
into the direct product of order-two groups:
$G/Q=\gp{d_1Q}_2\times\dots\times\gp{d_mQ}_2$ and consider
the words
$v_\chi(x_1,\dots,x_m,y)\in F(x_1,\dots,x_m,y)$
(in the free group) obtained from $w_\chi$
(see formula $(1')$) by substitution
$c_i$ with their
expressions via generators~$d_1,\dots,d_m$ and
substitution $q$ with a new letter $y$:
$$
v_\chi(d_1,\dots,d_m,q)=w_\chi(\~q).
$$
Formula (3) shows that the equation
$$
\prod_{\chi\in X}
\(v_\chi\(x_1,\dots,x_m,
\prod_{i=1}^n y_{\chi,i}^{2|T(Q)|}\)\)^{k_\chi}
=
\(a^2\)^{2^{|C|}\cdot|T(Q)|},
\eqno{(4)}
$$
where $n$ is a number such that each element of $Q$
is a product of $n$ squares (e.g., $n=\rk Q+1$),
has a solution in $G$:
$$
x_j=d_j,
\quad
y_{\chi,i}=g_{\chi,i},
\qbox{where $g_{\chi,i}\in G$ are such that
$\prod\limits_i g_{\chi,i}^2=q(\chi)$}.
$$
It remains to show that equation (4) has no solutions in
the dihedral group
$H=\gp b_2\semitimes \gp a_\infty$.

A substitution $x_j=b^{\epsilon_j}a^{k_j}$ naturally defines a
homomorphism $\phi\:C\to H/\gp{a^2}$ and a
character $\chi'\:C\to\{\pm1\}$:
$$
\chi'(d_jQ)=(-1)^{\epsilon_j},
\qbox{i.e. $\chi'=\^\chi\circ\phi$,
where $\^\chi\:H/\gp{a^2}\to\{\pm1\}$
is the
character of the action of $H/\gp{a^2}$
on $\gp{a^2}$.
}
$$

Let us substitute the
variables
$y_{\chi,i}$ by elements $h_{\chi,i}\in H$ and
note that
$$
\prod_{i=1}^n y_{\chi,i}^{2}=\prod_{i=1}^n h_{\chi,i}^2
=
a^{2l_\chi}
\qbox{for some $l_\chi\in\Z$}.
$$
Then the multiplicative
analogue of $(*)$ gives the equality
$$
v_\chi\(x_1,\dots,x_m,\prod_{i=1}^n y_{\chi,i}^{2|T(Q)|}\)
=
v_\chi\(x_1,\dots,x_m,a^{2l_\chi\cdot|T(Q)|}\)
=
\cases{\(a^{2l_{\chi'}}\)^{2^{|C|}\cdot|T(Q)|},
& if  $\chi=\chi'$;
\cr
\cr
1,
& if  $\chi\ne\chi'$.
}
$$
Therefore, this substitution transforms
the left-hand side of equation (4)
into
$$
a^{2l_{\chi'}\cdot2^{|C|}\cdot|T(Q)|\cdot k_\chi}
\ne
a^{2\cdot2^{|C|}\cdot|T(Q)|},
\qbox{because $k_\chi\ne\pm1$,}
$$
and (4) has no solutions in $H$ as required.

\s 4.
An example

Consider the following
example:
$$
G=D_{\infty}\times D_{\infty}=
\Bigl(\gp{b_1}_2\semitimes\gp{a_1}_\infty\Bigr)
\times
\Bigl(\gp{b_2}_2\semitimes\gp{a_2}_\infty\Bigr)
\hbox{ and }
G\supset H=\langle b \rangle_{2}\rightthreetimes
\langle a\rangle_{\infty}
\iso D_\infty,
\hbox{ where } b=b_{1}b_{2}\hbox{ and } a=a_{1}^{3}a_{2}^{5}.
$$
In this case,
$
Q=\gp{\{g^{2}\;|\;g\in G\}}=
\langle a_{1}^{2}\rangle_{\infty}
\times \langle a_{2}^{2}\rangle_{\infty}
\simeq
\Z\oplus\Z
$
and
$
C = G/Q =
\langle a_{1}Q \rangle_{2}
\times
\langle b_{1}Q \rangle_{2}
\times
\langle a_{2}Q \rangle_{2}
\times
\langle b_{2}Q \rangle_{2}
$.
Thus, there are $2^{4} = 16$ different characters
$C\to\{\pm1\}$,
and
only for two of them, $\alpha$ and $\beta$,
the subgroups~$Q_\chi$ are nontrivial:
$$
\alpha:b_{1}\mapsto-1,\
{a_{1}, a_{2}, b_{2}}\mapsto1,
\qquad
\beta:b_{2}\mapsto-1,\
{a_{1}, b_{1}, a_{2}}\mapsto1.
$$
Clearly, the element $a^2$ is not simple and we have to
prove that the subgroup $H$ is not verbally closed.

The lengthy word
$v_\chi(x_1,x_2,x_3,x_4,y)$
is the composition (in an arbitrary order) of
the following
16 words (as functions of $y$):
$$
\eqalign{
f_\chi(1,y),\
f_\chi(x_1,y),\
f_\chi(x_1x_2,y),\
f_\chi(x_1x_2x_3,y),\
f_\chi(x_1x_2x_3x_4,y),\
f_\chi(x_1x_2x_4,y),\
f_\chi(x_1x_3,y),\
f_\chi(x_1x_3x_4,y),
\cr
f_\chi(x_1x_4,y),\
f_\chi(x_2,y),\
f_\chi(x_2x_3,y),\
f_\chi(x_2x_3x_4,y),\
f_\chi(x_2x_4,y),\
f_\chi(x_3,y),\
f_\chi(x_3x_4,y),\
f_\chi(x_4,y).
}
\eqno{(**)}
$$
Here, the first arguments are
all expressions of the form
$x_1^{\epsilon_1}x_2^{\epsilon_2}x_3^{\epsilon_3}x_4^{\epsilon_4}$,
where $\epsilon_i\in\{0,1\}$,
and
$$
f_\chi(x_1^{\epsilon_1}x_2^{\epsilon_2}x_3^{\epsilon_3}x_4^{\epsilon_4},y)
=
y
\cdot
x_1^{\epsilon_1}x_2^{\epsilon_2}x_3^{\epsilon_3}x_4^{\epsilon_4}
\cdot
y^{\chi(a_1^{\epsilon_1}b_1^{\epsilon_2}a_2^{\epsilon_3}b_2^{\epsilon_4})}
\cdot
\(x_1^{\epsilon_1}x_2^{\epsilon_2}x_3^{\epsilon_3}x_4^{\epsilon_4}\)^{-1}.
$$
For example,
$
f_\alpha(x_1^{\epsilon_1}x_2^{\epsilon_2}x_3^{\epsilon_3}x_4^{\epsilon_4},y)
=
y
\cdot
x_1^{\epsilon_1}x_2^{\epsilon_2}x_3^{\epsilon_3}x_4^{\epsilon_4}
\cdot
y^{(-1)^{\epsilon_2}}
\cdot
\(x_1^{\epsilon_1}x_2^{\epsilon_2}x_3^{\epsilon_3}x_4^{\epsilon_4}\)^{-1}
$.
We see that, in the dihedral group,
$$
f_\alpha\Bigl(
\(a^{k_1}b^{\delta_1}\)^{\epsilon_1}
\(a^{k_2}b^{\delta_2}\)^{\epsilon_2}
\(a^{k_3}b^{\delta_3}\)^{\epsilon_3}
\(a^{k_4}b^{\delta_4}\)^{\epsilon_4}
,
a^{2k}\Bigr)
=
\cases{
a^{4k}, &if $\epsilon_2+\sum\limits_{i=1}^4\delta_i\epsilon_i$ is even;
\cr
\cr
1, &if $\epsilon_2+\sum\limits_{i=1}^4\delta_i\epsilon_i$ is odd;
}
\!\!\!\hbox{where $k_i\in\Z$ and $\delta_i,\epsilon_i\in\{0,1\}$.}
$$
Thus, if we take
$\delta_2=1$ and $\delta_1=\delta_3=\delta_4=0$,
then
$f_\alpha\Bigl(
\(a^{k_1}b^{\delta_1}\)^{\epsilon_1}
\(a^{k_2}b^{\delta_2}\)^{\epsilon_2}
\(a^{k_3}b^{\delta_3}\)^{\epsilon_3}
\(a^{k_4}b^{\delta_4}\)^{\epsilon_4}
,
a^{2k}\Bigr)
$
becomes equal to $a^{4k}$
for any choice of $\epsilon_i$;
if we take any other tuple of
$\delta_i\in\{0,1\}$, then at least one of 16
expressions $(**)$ becomes 1 after a substitution
$x_i\to a^{k_i}b^{\delta_i}$ and $y\to a^{2k}$.
This means that, for the composition
$v_\alpha$ of expressions $(**)$, we have
$$
v_\alpha\Bigl(
\(a^{k_1}b^{\delta_1}\)^{\epsilon_1},
\(a^{k_2}b^{\delta_2}\)^{\epsilon_2},
\(a^{k_3}b^{\delta_3}\)^{\epsilon_3},
\(a^{k_4}b^{\delta_4}\)^{\epsilon_4},
a^{2k}
\Bigr)
=
\cases{
a^{2^{17}k},
& if $(\delta_1,\delta_2,\delta_3,\delta_4)=(0,1,0,0)$;
\cr
1,
& otherwise.
}
$$
All other characters
behave
similarly.
For instance,
$$
v_{\alpha\beta}\Bigl(
\(a^{k_1}b^{\delta_1}\)^{\epsilon_1},
\(a^{k_2}b^{\delta_2}\)^{\epsilon_2},
\(a^{k_3}b^{\delta_3}\)^{\epsilon_3},
\(a^{k_4}b^{\delta_4}\)^{\epsilon_4},
a^{2k}
\Bigr)
=
\cases{
a^{2^{17}k},
& if $(\delta_1,\delta_2,\delta_3,\delta_4)=(0,1,0,1)$;
\cr
1,
& otherwise.
}
$$
Equations (4) has the form
$$
\Bigl(v_\alpha\(x_1,x_2,x_3,x_4,
y_{\alpha,1}^2\)\Bigr)^3
\cdot
\Bigl(v_\beta\(x_1,x_2,x_3,x_4,
y_{\beta,1}^2\)\Bigr)^5
\cdot
\prod_{\chi\ne\alpha,\beta}
\Bigl(v_\chi\(x_1,x_2,x_3,x_4,
y_{\chi,1}^2\)\Bigr)^{\the\year}
=
a^{2^{17}}.
$$
(We take the exponent \the\year\ to emphasize that
we can put any number here, except $\pm1$; of course, the simplest choice
is to
replace \the\year\ with 0.)

What is said means that, for any substitution
$x_i\to a^{k_i}b^{\delta_i}$, the
left-hand side of the equation
takes a value in
$$
\eqalign{
&\gp{a^{2^{17}\cdot3}},
\hbox{ if $(\delta_1,\delta_2,\delta_3,\delta_4)=(0,1,0,0)$};
\cr
&\gp{a^{2^{17}\cdot5}},
\hbox{ if $(\delta_1,\delta_2,\delta_3,\delta_4)=(0,0,0,1)$};
\cr
&\gp{a^{2^{17}\cdot\the\year}}
\hbox{ in all other cases}.
}
$$
Thus, equation (4) has no solutions in the dihedral group.
On the other hand, in $G$, there is a solution:
$x_1=a_1$,
$x_2=b_1$,
$x_3=a_2$,
$x_4=b_2$,
$y_{\alpha,1}=a_1$,
$y_{\beta,1}=a_2$,
$y_{\chi,1}=1$ for $\chi\notin\{\alpha,\beta\}$.


\REFERENCES

\def\textbf#1{{\bf#1}}

\[KlMa18]
Ant. A. Klyachko, A. M. Mazhuga,
Verbally closed virtually free subgroups,
\emph{Sbornik: Mathematics} {\bf 209} (2018) (to appear).
\newline
See also arXiv:1702.07761.

\[RKh13]
V. A. Roman'kov, N. G. Khisamiev,
Verbally and existentially closed subgroups of free nilpotent groups.
\emph{Algebra and Logic} \textbf{52} (2013), 336-351.

\[Lee02]
D. Lee,
On certain C-test words for free groups.
\emph{J. Algebra}, \textbf{247} (2002), 509-540.

\[Mazh17]
A. M. Mazhuga,
On free decompositions of verbally closed subgroups
of free products of finite groups,
\emph{J. Group Theory}, {\bf 20:5} (2017), 971-986.
\newline
See also arXiv:1605.01766.

\[Mazh18]
A. M. Mazhuga,
Strongly verbally closed groups,
\emph{J. Algebra}, {\bf 493} (2018), 171-184.
\newline
See also
arXiv:1707.02464

\[MR14]
A. Myasnikov, V. Roman'kov,
Verbally closed subgroups of free groups,
\emph{J. Group Theory}, \textbf{17} (2014), 29-40.
\newline
See also arXiv:1201.0497.

\[Rom12]
V. A. Roman'kov,
Equations over groups,
\emph{Groups Complexity Cryptology},
\textbf{4:2} (2012), 191-239.

\[Sta71]
J. Stallings, Group theory and three-dimensional manifolds,
Yale Math. Monographs (1971).

\end